\input amssym.def
\input amssym
\magnification=1200
\parindent0pt
\hsize=16 true cm
\baselineskip=13  pt plus .2pt
$ $

\def\Z{\Bbb Z}
\def\D{\Bbb D}
\def\A{\Bbb A}

\centerline {\bf On finite simple groups acting on homology spheres}

\centerline {\bf with small fixed point sets}

\bigskip

\centerline {Bruno P. Zimmermann}

\bigskip

\centerline {Universit\`a degli Studi di Trieste}

\centerline {Dipartimento di Matematica e Geoscienze}

\centerline {34127 Trieste, Italy}

\bigskip \bigskip

Abstract.  {\sl A finite nonabelian simple group does not admit a free action
on a homology sphere, and the only finite simple group which acts on a homology
sphere with at most 0-dimensional fixed point sets ("pseudofree action") is the
alternating group $\A_5$ acting on the 2-sphere. Our first main theorem is the
finiteness result that there are only finitely many finite simple groups which
admit a smooth action on a homology sphere with at most $d$-dimensional fixed
points sets, for a fixed $d$.  We then go on proving that  the finite simple
groups acting on a homology sphere with at most 1-dimensional fixed point sets
are the alternating group $\A_5$ in dimensions 2, 3 and 5,  the linear
fractional group ${\rm PSL}_2(7)$ in dimension 5, and possibly the unitary
group ${\rm PSU}_3(3)$ in dimension 5 (we conjecture that it does not admit any
action on a homology 5-sphere but cannot exclude it at present). Finally, we
discuss the situation for arbitrary finite groups which admit an action on a
homology 3-sphere. }

\bigskip \bigskip

{\bf 1. Introduction}

\medskip

We are interested in finite groups, and in particular in finite simple groups, which admit
a smooth action on an integer or a mod 2 homology sphere. A homology sphere (resp. a
mod 2 homology sphere) is a closed manifold  with the  integer homology of a
sphere (resp. the  mod 2 homology of a sphere, i.e. with coefficients in the integers mod
2).  In the present paper, simple group will always mean {\it nonabelian} simple group;
also, all actions will be smooth (or locally linear), orientation-preserving and faithful.

\medskip

By  [24], [15-17], the only finite simple group which admits an action on a homology
3-sphere is the alternating group $\A_5$, and the only finite simple groups acting on a
homology 4-sphere are the alternating groups $\A_5$ and $\A_6$. The finite simple groups
acting on a homology 5-sphere are
considered in [10, Theorem 2] (see Theorem 4 in section 2).

\medskip

Suppose that a finite nontrivial group $G$ admits a {\it free}  action on an integer
homology sphere of dimension $n$.  Since we are considering orientation-preserving actions,
by the Lefschetz fixed point theorem this is possible only in odd dimensions; also,
$G$ has periodic cohomology (cf. [4, chapters I.6 and VII.10]), and the class of groups of
periodic cohomology is well-known and quite restricted (see [1]); in particular, no finite
simple groups occur.  A next case which has been considered is that of {\it pseudo-free
actions}, i.e.  actions with 0-dimensional fixed point sets; such actions exist only
in even dimensions since the fixed point set of any finite cyclic (orientation-preserving)
subgroup $\Bbb Z_p$ has even codimension (since this is the case for the linear action
induced on the tangent space of a fixed point). For this case it is easy to see that
again only the groups with periodic cohomology occur, plus the finite groups acting on the
2-sphere (see the Remark at the end of section 2); in particular, the only finite simple
group which occurs is the alternating group $\A_5$ acting on $S^2$.

\medskip

Thus one is led to consider less restrictive conditions on fixed point sets.
For an integer $d \ge -1$, we say that a finite group acts {\it  with at most
$d$-dimensional fixed points sets} if the fixed point of each nontrivial
element has dimension at most $d$ (where $d = -1$ stands for empty fixed point
set). Our first main  result is the following:

\bigskip

{\bf Theorem 1.}  {\sl For a fixed $d$, there are only finitely many finite simple
groups which admit an action  on a homology sphere with at most $d$-dimensional fixed
point sets.}

\bigskip

We believe that such a finiteness result does not hold for actions of finite simple
groups on mod 2 homology spheres. For example, it is likely that all groups
${\rm PSL}_2(p)$, $p$ prime, admit an action already on a mod 2 homology 3-sphere, that is
in dimension three; examples of such actions for various small values of $p$ are given in
[26] but the problem remains open in general (see also the survey [25]).

\medskip

Next we consider the case $d=1$:

\bigskip

{\bf Theorem 2.}  {\sl  The finite simple groups which act with at most 1-dimensional
fixed point sets  on a homology sphere are the alternating group $\A_5$ in
dimensions  2, 3 and 5, the linear fractional group ${\rm PSL}_2(7)$ in dimension
5, and possibly the unitary group ${\rm PSU}_3(3)$ in dimension 5.}

\bigskip

We note that the groups $\A_5$ and ${\rm PSL}_2(7)$ admit linear actions with
at most 1-dimensional fixed point sets on spheres of the indicated dimensions
(see the proof of Theorem 2). The unitary group ${\rm PSU}_3(3)$ has a linear
action on $S^6$ (with at most 2-dimensional fixed point sets); we conjecture
that it does not admit any action on a homology 5-sphere but cannot exclude it
at present.

\medskip

The proofs of Theorems 1 and 2 are based also on a part of the following
theorem which collects some consequences of the Borel-formula for actions of an
elementary abelian $p$-group ([2, Theorem XIII.2.3]), in combination with some  deep
results from the theory of finite simple groups.

\bigskip

{\bf Theorem 3.} {\sl i) Let $G$ be a finite group which admits an action on a
mod 2 homology $n$-sphere such that involutions have at most $d$-dimensional fixed
points sets. Suppose that $G$ has a subgroup isomorphic to the Klein group $\Bbb Z_2
\times \Bbb Z_2$. Then  $$\; n \le 3d+2;$$ in particular, this holds if $G$ is a nonabelian
simple group.

\medskip

ii) Suppose that $n = 3d+2$. If $G$ is a nonabelian simple group then $G$ has
2-rank at most two (i.e., $G$ has no subgroups $(\Z_2)^3$) and is one of the following
groups (where $q$ denotes an odd prime power):
$${\rm PSL}_2(q), \;\; {\rm PSL}_3(q), \;\; {\rm PSU}_3(q), \;\;
{\rm PSU}_3(4), \;\; \A_7, \;\; {\rm the \; Mathieu \; group} \; M_{11}.$$

iii)  Suppose that $n > 3d+2$. Then $G$ has no subgroup $(\Z_2)^2$,  and  hence
2-periodic cohomology. Moreover, if $G$ is nonsolvable then it has the following
structure: denoting by ${\cal O}(G)$ the maximal normal subgroup of odd order
of $G$, the factor group $G/{\cal O}(G)$ contains a normal subgroup  of odd
index which is isomorphic to
$${\rm SL}_2(q), \;\; {\rm TL}_2(q) \;\; or  \;\; \hat\Bbb A_7.$$}

Here $\hat\Bbb A_7$ denotes the unique perfect central extension of the alternating group
$\A_7$, with center of order two, and ${\rm TL}_2(q)$
is the 2-fold extension of ${\rm SL}_2(q)$ with a unique involution  whose quotient group
is isomorphic to ${\rm PGL}_2(q)$ (see [1, chapter IV.6]).

\medskip

We refer to [4, Theorem VI.9.7] and [23] for the notion of a $p$-periodic group.
We note also that a finite group admits a free action on a mod 2 homology sphere if and
only if it is 2-periodic and has a unique involution, see [20].

\medskip

In section 3 we consider finite groups acting on a homology 3-sphere.  We note
that, by the recent geometrization of finite group-actions on 3-manifolds due
to Thurston and Perelman, the finite groups which admit an action on the
3-sphere $S^3$ are exactly the finite subgroups of the orthogonal group O(4).
This is no longer true for finite groups which admit an action on an arbitrary
homology 3-sphere; however the classification of such groups remains open and
appears to be difficult (even for quite easy types of finite groups as in
Question 2 i) of section 3).

\bigskip

{\bf 2. Proofs}

\medskip

We start with the {\it Proof of Theorem 3}. Let $G$ be a finite group acting on a mod
2 homology $n$-sphere with at most $d$-dimensional fixed point sets; suppose that $G$ has a
subgroup $(\Z_2)^2$. We note that in particular every finite nonabelian simple group has a
subgroup $(\Z_2)^2$: in fact, if a finite group has no subgroup
$(\Z_2)^2$  then, by a theorem of Burnside ([22, 4.4.3] or [4, Theorem IV.4.3]), a Sylow
2-subgroup is either cyclic or generalized quaternion, but by [22, p. 144,
Corollary 1] and [22, p. 306, Example 3] a Sylow 2-subgroup of a finite simple
group cannot be cyclic or generalized quaternion.

\medskip

By Smith fixed point theory ([3], [21]), the fixed point set of an
orientation-preserving periodic map of prime order $p$ of a mod $p$ homology
sphere is again a mod $p$ homology sphere, of even codimension.
A basic tool for actions of an elementary abelian $p$-group $A \cong (\Z_p)^m$ on a mod
$p$  homology $n$-sphere is then the following Borel formula ([2, Theorem XIII.2.3]):

$$n-r =  \sum _{H} \; (r(H)-r)$$

where the sum is taken over all subgroups
$H$ of index $p$ of $A$, $r(H)$ denotes the dimension of the fixed point set of a given
subgroup $H$ and $r = r(A)$ the dimension of the fixed point set of $A$ (equal to -1
if the fixed point set is empty).

\medskip

Applying the Borel formula to a subgroup $A \cong (\Z_2)^2$ of $G$ and using the fact that
$-1 \le r(A), r(H) \le d$, we obtain the inequality
$$n = \sum _{H} r(H) \, - \, 2r \le  3d + 2,$$
proving part i) of Theorem 3.

\medskip

Suppose that $n = 3d + 2$. Then $r(H) = d$ for each of the three $\Z_2$ subgroups of $A$,
and $r = r(A) = -1$.  Suppose that $G$ has a subgroup
$B \cong (\Z_2)^3$. Then $B$ has exactly seven subgroups $A \cong (\Z_2)^2$ of
index two, $r(A) = -1$ for each of these by the above, and in particular also
$r(B) = -1$. Applying the Borel formula to $B$ now  we obtain a contradiction ($n = -7
-6r = -13$).  So
$G$ does not have an elementary abelian subgroup $(\Z_2)^3$ of rank three and has
2-rank equal to two. By  a fundamental result in the classification of the
finite simple groups ([9]), the finite simple groups of 2-rank two are exactly the groups
listed in  Theorem 3 ii).

\medskip

Finally, suppose that $n > 3d + 2$. Then as before  $G$ has no
subgroups  $(\Z_2)^2$ and hence, by [4, Theorem VI.9.7], $G$ has 2-periodic
cohomology.

\medskip

Suppose that $G$ is nonsolvable. By the Feit-Thompson theorem, a Sylow
2-subgroup $S$ of $G$ in nontrivial. Since the finite 2-group $S$ has
nontrivial center and no subgroups $(\Z_2)^2$, $S$ has a unique involution. By
the theorem of Burnside above, a finite 2-group with a unique involution is
either cyclic or generalized quaternion. Since the Sylow 2-subgroup of a
nonsolvable group cannot be cyclic ([22, chapter 5.2, Corollay 2]), $S$ is a
generalized quaternion group. The  structure of the finite nonsolvable groups
with a generalized quaternion Sylow 2-subgroup is given in  [22, chapter 6,
Theorem 8.7], and the version given in Theorem 3 is an elaboration of this as
in [25, Theorem 5].

\medskip

This completes the proof of Theorem 3.

\bigskip

Theorem 1 is now a consequence of Theorem 1 i) and of [10, Theorem 1] stating
that for each dimension $n$ there are only finitely many finite
simple groups which admit an action on a homology $n$-sphere.

\bigskip

For the {\it Proof of Theorem 2}, let $G$ be a finite simple group acting on a homology
$n$-sphere with at most 1-dimensional fixed point sets; by Theorem 3 i),  $n \le 5$.

\medskip

By [24], the only finite simple group acting on a homology 3-sphere is $\A_5$.

\medskip

As noted above, by Smith fixed point theory the fixed point set of an
orientation-preserving periodic map of prime order $p$ of a mod $p$ homology
sphere is again a mod $p$ homology sphere of even codimension. Now the case $n
= 4$ is excluded by the Borel formula applied to a subgroup $(\Z_2)^2$ of $G$
(since the hypothesis of at most 1-dimensional fixed point sets implies that
$r(\Z_2) = 0$ in dimension 4).

\medskip

So we are left with the case $n = 5$: suppose that $G$ acts on a homology
5-sphere.  By Theorem 3 ii), $G$ has 2-rank two and is one of the following
groups:
$${\rm PSL}_2(q), \;\; {\rm PSL}_3(q), \;\; {\rm PSU}_3(q), \;\; {\rm
PSU}_3(4), \;\; \A_7, \;\; {\rm the \; Mathieu \; group} \; M_{11}.$$

For the proof of Theorem 2 we have to exclude all of these groups except $\A_5
\cong {\rm PSL}_2(5)$, ${\rm PSL}_2(7)$ and ${\rm PSU}_3(q)$. This is based on
the following result from [10] (resp. on some part of its proof):

\bigskip

{\bf Theorem 4.} ([10, Theorem 2]) {\sl  The finite simple groups which admit an
action on a homology 5-sphere are the following:

$$\A_5 \cong  {\rm PSL}_2(5), \;\; \A_6 \cong {\rm PSL}_2(9), \;\;\A_7, \;\; {\rm
PSL}_2(7), \;\; {\rm PSU}_4(2), \;\; possibly \;\; {\rm PSU}_3(3).$$}

With the exception of the unitary group ${\rm PSU}_3(3)$ (which admits a linear
action on $S^6$) these are exactly the finite
simple groups which admit a linear action on $S^5$. We note that the proof of [10, Theorem
2] is on the basis of the full classification of the finite simple groups; in our situation
the proof is much shorter since we have to consider only the quite restricted list of the
finite simple groups of 2-rank two.

\medskip

Note that for the proof of Theorem 2 we still have to exclude that the
alternating groups $\A_6$ and $\A_7$ from the list in Theorem 4 admit an action on a
homology 5-sphere with at most 1-dimensional fixed point sets. Suppose that $\A_6$ admits
such an action. We consider an elementary abelian subgroup $(\Z_3)^2$ of $\A_6$ generated by
two disjoint cycles of length three, with  four subgroups $\Z_3$; note that the four
subgroups
$\Z_3$ of $(\Z_3)^2$ are conjugate in pairs in $\A_6$ (two are generated by a 3-cycle, the
other two by a product of two 3-cycles).  By the Borel formula, $n-r = 5-r = \sum _{i} \;
(r(H_i)-r)$, or

$$5+3r = r(H_1) +  r(H_2) +  r(H_3) +  r(H_4)$$

where the $H_i$ denote the four
subgroups $\Z_3$ of $(\Z_3)^2$. By our assumption of at most 1-dimensional fixed
point sets, we have that $r(H_i) \le 1$, hence also $r \le 1$ and $5+3r \le 4$. This
excludes  the possibilities $r=0$ and $r=1$. Suppose that $r = -1$. Then the only solution
of the Borel formula is
$5+3r = 2 = r(H_1) +  r(H_2) +  r(H_3) +  r(H_4) =  1+1+1-1$ (note that $r(H_i) \ne 0$
since the fixed point set of each $H_i$ has even codimension); however also this solution
is not possible since the four subgroups
$H_i$ of $(\Z_3)^2$ are conjugate in pairs.

\bigskip

This excludes the groups $\A_6$ and $\A_7$, and
the only finite simple groups which remain from the list in Theorem 4
are $\A_5$, ${\rm PSL}_2(7)$ and
${\rm PSU}_3(3)$. The dodecahedral group $\A_5$ acts on $S^2$ with
0-dimensional fixed point sets $S^0$, and it has two linear actions on $S^3$:
one is the suspension of the action on $S^2$, with two global fixed points, the
other one is the restriction to $S^3$ of its irreducible 4-dimensional real
representation  (a summand of the standard 5-dimensional representation of
$\A_5$ by permutation of coordinates of $\Bbb R^5$). A linear action of $\A_5$
on $S^5$ is obtained by considering $S^5\cong S^2 * S^2$ as the join of two
2-spheres and by taking also the join of two actions of $\A_5$ on $S^2$, with
fixed point sets $S^0 * S^0 \cong S^1$ (or equivalently by restricting to $S^5$
the direct sum of two irreducible 3-dimensional real representations).  The
group ${\rm PSL}_2(7)$ has an irreducible  3-dimensional complex
representation, and the restriction of the corresponding 6-dimensional
(reducible) real representation to $S^5$ has 1-dimensional fixed point sets
(see the character tables in [5]; note that ${\rm PSL}_2(7)$ has also an
irreducible 6-dimensional real representation whose restriction, however, has
also 3-dimensional fixed point sets).

\medskip

This completes the proof of Theorem 2.

\bigskip

{\bf Remark.}  Suppose that a finite nontrivial group $G$ acts
orientation-preservingly and pseudofreely (i.e., with at most $0$-dimensional fixed
point sets) on a homology $n$-sphere.  Applying the Borel formula,
$G$ has no subgroup $(\Z_2)^2$ if $n > 2$, and no subgroup $(\Z_p)^2$ for odd primes $p$.
So, if $n > 2$, every abelian subgroup of $G$ is cyclic and hence $G$ has periodic
cohomology (see [4, Proposition VI.9.3]). The groups of periodic cohomology are
well-known (see e.g. [1]). If such a group has in addition a unique involution then it is
known as an application of  high-dimensional  surgery theory that it admits a free action
on a sphere of odd dimension ([14]); by suspending such a free action one obtains a
pseudofree action with exactly two global fixed points on a sphere of even dimension. Thus
the finite groups which admit a pseudofree action on some homology sphere are exactly the
finite  groups acting on $S^2$, plus the groups of periodic cohomology with a unique
involution. Kulkarni has shown ([12, Theorem 7.4]) that, with the only exception of maybe
dihedral groups, every pseudofree action on a homology sphere of dimension at least three
has exactly two global fixed points such that the action on the complement of these two
fixed points is free ("semifree action").

\bigskip
\vfill \eject

{\bf 3.  The situation in dimension three}

\medskip

By the  recent geometrization of finite group-actions on 3-manifolds due to
Thurston and Perelman, every finite group of diffeomorphisms of $S^3$ is
conjugate to a subgroup of the orthogonal group O(4); in particular, the finite
groups occurring are exactly the finite subgroups of the orthogonal group O(4).
The finite groups which admit an action on an arbitrary homology 3-sphere are
discussed in [27]; a complete classification of these groups is not known.  We consider
first the case of free actions.

\medskip

If a finite group $G$ admits a {\it free} action on a homology 3-sphere then $G$ has
periodic cohomology of period four and a unique involution; a list of such groups is
given in [19], together with the subclass of all finite groups which admit a free,
{\it linear} action on the 3-sphere. By [13] there remains one class of groups
$Q(8a,b,c)$ in [19] which do not admit a free, linear action on $S^3$ but for which
the existence of a free action on a homology 3-sphere remains open, in general.
By [18], some of the groups  $Q(8a,b,c)$ admit a free
action on a homology 3-sphere and some others do not, but the exact classification remains
open  (see also the discussion in [11, Problem 3.37 Update A (p.173)]).

\medskip

The group $Q(8a,b,c)$ has a presentation
$$< x,y,z \,|\, x^2=(xy)^2=y^{2a}, z^{bc}=1, xzx^{-1}=z^r, yzy^{-1}=z^{-1} >, $$

for relatively coprime positive integers $8a, b$ and $c$ such that either $a$
is odd and $a>b>c\ge 1$, or
$a \ge 2$ is even and $b>c\ge 1$; also, $r \equiv -1$ mod $b$ and $r
\equiv +1$ mod $c$. Note that
$Q(8a,b,c)$ is a semidirect product $\Z_{bc} \rtimes Q(8a)$, with normal subgroup
$\Z_{bc} \cong \Z_b \times \Z_c$ generated by $z$, and factor group the
generalized quaternion or binary dihedral subgroup $Q(8a) \cong Q(8a,1,1)$ of
order $8a$ generated by $x$ and $y$. See also [6, section 7] for a description of these
groups and various inclusions between them.

\medskip

We note that a group
$Q(8a,b,c)$ does not admit a free action on a homology 3-sphere if $a$ is even ([13], [6]).
If $a$ is odd then $Q(8a) \cong \Z_a \rtimes Q(8)$, and hence
$Q(8a,b,c)$ is an extension of $\Z_{abc}  \cong \Z_a \times \Z_b \times \Z_c$ by the
quaternion group $Q(8) = \{\pm 1,\pm i,\pm j,\pm k\}$ such that  $i, j, k$ acts
trivially on $\Z_a$, $\Z_b$, $\Z_c$,  respectively, and in a dihedral way on the other two.

\medskip

Concerning {\it nonfree} actions of the groups $Q(8a,b,c)$, we note the following:

\bigskip

{\bf Proposition 1.}  {\sl A group $Q(8a,b,c)$ does not admit a  nonfree action on a
homology 3-sphere (orientation-preserving or not).}

\medskip

{\it Proof.}  Suppose that $G = Q(8a,b,c)$ acts orientation-preservingly on a
homology 3-sphere $M$. The unique involution
$h = x^2=(xy)^2=y^{2a}$ of $G$ is central in $G$; by Smith fixed point theory
the fixed point set of $h$ has even codimension and is either empty or a 1-sphere $S^1$ in
$M$ (see e.g. [3]). Suppose that $h$ has nonempty fixed point set
$S^1$; note that $S^1$ is invariant under the action of $G$. We note that, if a finite
orientation-preserving group leaves invariant a 1-sphere $S^1$ in a 3-manifold then $G$ is
isomorphic to a subgroup of a semidirect product
$(\Z_m \times \Z_n) \rtimes \Z_2$, with a dihedral action of $\Z_2$ on
$\Z_m \times \Z_n$ (here $\Z_2$ acts as a reflection or strong inversion on
$S^1$ whereas $\Z_m \times \Z_n$ acts by rotations about and along $S^1$). Since clearly
$G = Q(8a,b,c)$ is not of this type, the unique involution
$h$ of $G$ has to act freely, and hence also every nontrivial element in $G$ of
even order.

\medskip

Next suppose that some nontrivial element $g \in G$ of odd prime order  has
nonempty fixed point set $S^1$, acting as a rotation about $S^1$;  we
can assume that $g$ is an element of one of the subgroups $\Z_b \times \Z_c$ or $Q(8a) \cong
Q(8a,1,1)$ of $G$. In each case some element $u$ of even order
in the generalized quaternion group $Q(8a)$ acts dihedrally on $g$ (i.e.,
$ugu^{-1}=g^{-1}$).  Since $u$ has no fixed points, it acts
as a rotation along and about $S^1$. But then $u$ commutes with the rotation $g$
about $S^1$ which is a contradiction.

\medskip

Now suppose that some element of $G$ reverses the orientation of $M$. Since the order $bc$
of $z$ is odd, $x$ or $xy$ are orientation-reversing; we assume
that $x$ is orientation-reversing (the case of $xy$ is analogous).  By Smith fixed point
theory, the fixed point set of $x$ has odd codimension and is either a 0-sphere (two
points) or a 2-sphere. Then the fixed point set of the central involution $h = x^2$ is also
nonempty and hence a 1-sphere $S^1$ (of even codimension), the fixed point set
of $x$ is a 0-sphere $S^0 \subset S^1$, and $x$ acts as a reflection (strong inversion) on
$S^1$.

\medskip

If both $x$ and $xy$ reverse the orientation of $M$ then also $xy$ acts as a
reflection on $S^1$, and $y$ is orientation-preserving and acts as a rotation about and
along $S^1$. But also the subgroup $\Z_{bc}$ generated by $z$ acts as a group of rotations
about and along $S^1$, so $y$ and $z$ commute; this is a contradiction since $y$ acts
dihedrally on $z$.

\medskip

So $xy$ acts orientation-preservingly and its fixed point set is either empty
or $S^1$. If $xy$ has empty fixed point set then it acts by rotations about and
along $S^1$, and hence commutes with $y^2$ (of order $2a \ge 4$) which acts
also by rotations about and along $S^1$. This is a contradiction since $x$ and
hence $xy$ acts dihedrally on $y$ and $y^2$ (specifically, $(xy)^2 = xyxy = h$
implies that $xyx^{-1} = hy^{-1}x^{-2} = hy^{-1}h^{-1} = y^{-1}$).  If $xy$
fixes $S^1$ instead then it acts by rotations about $S^1$ and commmutes again
with $y^2$, so we obtain the same contradiction as before.

\medskip

So there are no orientation-reversing actions of $G$ on a homology 3-sphere.
This completes the proof of Proposition 1.

\bigskip

By results of Milgram [18], some of the groups $Q(8a,b,c)$ admit a free action
on a homology 3-sphere. By the geometrization of 3-manifolds with finite
fundamental group, none of the groups $Q(8a,b,c)$ admits a free action on $S^3$
(since, by [19], they do not admit a free, linear action on $S^3$). By
Proposition 1, they also don't admit nonfree actions on $S^3$ (alternatively,
considering orthogonal actions, one can confront them with the list of the
finite subgroups of O(4) in [8], see also [7] for the geometry of their
quotient orbifolds in the orientation-preserving case). Summarizing, the
following holds:

\bigskip

{\bf Proposition 2.}  {\sl The class of finite groups which admit an action
on a homology 3-sphere is strictly larger than the class of finite groups which admit an
action on $S^3$ (or the class of finite subgroups of ${\rm SO}(4)$).}

\bigskip

There arises naturally the question of how big the difference is between the classes of
groups in Proposition 2: do there occur other groups than the Milnor groups $Q(8a,b,c)$?
If a finite group admits a free action on a homology 3-sphere but not on $S^3$ then it is
in fact one of the Milnor groups $Q(8a,b,c)$, with $a$ odd ([13]), so any other such group
would admit only nonfree actions on a homology 3-sphere.

\bigskip

{\bf Question 1.}  i) Does there exist a finite group which admits a nonfree,
orientation-preserving  action on a homology 3-sphere but is not
isomorphic to a subgroup of the orthogonal group SO(4)?

\medskip

ii) Is there a finite group with an orientation-reversing action on a homology
3-sphere which is not isomorphic to a subgroup of O(4)?

\bigskip

In the following, concentrating on the orientation-preserving case, we discuss some natural
candidates. It is shown in [27] that the finite nonsolvable groups which admit an
orientation-preserving action on a homology 3-sphere are exactly the finite nonsolvable
subgroups of the orthogonal group SO(4) $\cong S^3 \times _{\Z_2} S^3$ (the central product
of two copies of the unit quaternions), plus possibly two other classes of groups:

\medskip

-  the central products
$$\A_5^* \times _{\Z_2}Q(8a,b,c)$$ where $a$ is odd and $\A_5^*$ denotes the binary
dodecahedral group;

\medskip

- their subgroups
$$\A_5^* \times _{\Z_2} (\D_{4a}^* \times \Z_b)$$

where $\D_{4a}^* \cong Q(4a) \cong \Z_a \rtimes \Z_4$ denotes the binary dihedral or
generalized quaternion group of order $4a$.

\bigskip

In turn these have subgroups
$$\D_8^* \times _{\Z_2} (\D_{4a}^* \times \Z_b)$$
which do not act freely on a homology 3-sphere (since they have a subgroup $\Z_2 \times
\Z_2$).

\bigskip

{\bf Lemma.}  {\sl For odd coprime integers $a,b \ge 3$, the group
$G =  \D_8^* \times _{\Z_2} \D_{4a}^* \times \Z_b$ does not admit an
orientation-preserving, linear action on $S^3$.}

\bigskip

{\it Proof.} Suppose that  $G$ is a subgroup of the orthogonal group SO(4)
$\cong S^3 \times_{\Z_2} S^3$. The finite subgroups of the unit quaternions
$S^3$   are cyclic, binary dihedral or binary polyhedral groups. The two
projections of the subgroup $\D_8^*$ of $G$ to the first and second factor
of $S^3 \times_{\Z_2} S^3$ are cyclic or binary dihedral groups; since $\D_8^*$
is nonabelian, one of the two projections, say the first one,  has to be a
binary dihedral group. Then, since the projections of the subgroups $\D_8^*$ and
$\D_{4a}^*$ of $G$ commute elementwise, the
projection of  $\D_{4a}^*$ to the second factor of $S^3 \times_{\Z_2} S^3$ has
to be a binary dihedral group.  But then at least one  of the two projections of the cyclic
subgroup $\Z_b$ of $G$ (any nontrivial one) does not commute elementwise with either
the binary dihedral projection of
$\D_8^*$ or that of $\D_{4a}^*$. This contradiction completes the proof of the Lemma.

\bigskip

{\bf Question 2.}  i) For odd, coprime integers  $a,b > 1$, does
$$\D_8^* \times_{\Z_2} \D_{4a}^* \times \Z_b$$
admit an orientation-preserving action on a homology 3-sphere? (If $a$
is even then there is no such action by [27, Lemma].)

\medskip

ii) Does the central product $$\Z_4  \times _{\Z_2} Q(8a,b,c)$$ admit an action on some
homology 3-sphere (assuming that $Q(8a,b,c)$ does)?

\bigskip

Note that these groups do not act freely on a homology
3-sphere (since they have a subgroup $\Z_2 \times \Z_2$) and are not
isomorphic to a subgroup of SO(4)  (by the Lemma, and since $Q(8a,b,c)$ is not).

\medskip

If the answer to i) is negative then by [27] the class of the finite
nonsolvable groups which admit an orientation-preserving action on a homology
3-sphere coincides with the class of the finite  nonsolvable subgroups of the
orthogonal group SO(4). On the other hand, if a group in i) or ii)  admits such
an action then this would give a first example of a finite group which admits a
nonfree, orientation-preserving action on a homology 3-sphere but which is not
isomorphic to a subgroup of SO(4) (and in case i) independently of the quite
difficult Milnor groups $Q(8a,b,c)$).

\medskip

Finally, considering also the case of mod 2 homology 3-spheres, we close with the following:

\bigskip

{\bf Conjecture.}  Each linear fractional group ${\rm PSL}_2(p)$, $p$ prime, admits an
action on a mod 2 homology 3-sphere.

\bigskip

By [15] or[16], these are exactly the candidates among the finite nonabelian simple groups
which possibly admit an action on a mod 2 homology 3-sphere.  Examples of such actions for
various small values of $p$ are given in [26]; see [15],[16] or the survey [25] for a partial
classification of the finite nonsolvable groups which admit an action on a mod 2 homology 3-sphere.

\bigskip

{\bf Acknowledgment.}  The author was supported by a FRV grant from
Universit\`{a} degli Studi di Trieste.

\bigskip

\centerline {\bf References}

\bigskip

\item {[1]}  A. Adem, R.J. Milgram, {\it Cohomology of finite groups.} Grundlehren der
Math. Wissenschaften 309, Springer-Verlag 1994

\item {[2]} A. Borel, {\it Seminar on Transformation Groups.} Ann. Math.  Studies 46,
Princeton University Press 1960

\item {[3]} G. Bredon, {\it Introduction to Compact Transformation Groups.} Academic
Press, New York 1972

\item {[4]} K.S. Brown, {\it Cohomology of Groups.}  Graduate Texts in Mathematics
87, Springer 1982

\item {[5]} J.H. Conway, R.T. Curtis, S.P. Norton, R.A. Parker, R.A.Wilson, {\it
Atlas of Finite Groups.} Oxford University Press 1985

\item {[6]} J.F. Davis, {\it The surgery semicharacteristic.}  Proc. London Math. Soc. 47,
411-428  (1983)

\item {[7]} W.D. Dunbar, {\it Geometric orbifolds.}  Rev. Mat. Univ. Complut. Madrid 1,
67-99  (1988)

\item {[8]}  P. Du Val, {\it  Homographies, Quaternions and Rotations.} Oxford
Math. Monographs, Oxford University Press 1964

\item {[9]} D. Gorenstein, {\it The Classification of Finite Simple Groups.} Plenum
Press, New York 1983

\item {[10]} A. Guazzi, B. Zimmermann, {\it On finite simple groups acting on
homology spheres.}  Mo-natsh. Math. 169,  371-381 (2013)

\item {[11]} R. Kirby, {\it Problems in low-dimensional topology.}  Geometric
Topology. AMS/IP Studies in Advanced Mathematics Volume 2, part 2, 35-358
(1997)

\item {[12]} R.S. Kulkarni, {\it Pseudofree actions and Hurwitz's $84(g - 1)$
theorem.}  Math. Ann. 261, 209 - 226  (1982)

\item {[13]} R. Lee, {\it Semicharacteristic classes.}   Topology 12, 183-199
(1973)

\item {[14]} I. Madsen, C.B. Thomas, C.T.C. Wall, {\it The topological
space form problem II: Existence of free actions.} Topology 15, 375-382
(1976)

\item {[15]} M. Mecchia, B. Zimmermann, {\it On finite groups acting on
$\Z_2$-homology 3-spheres.} Math. Z. 248, 675-693 (2004)

\item {[16]} M. Mecchia, B. Zimmermann, {\it On finite simple groups acting on integer
and mod 2 homology 3-spheres.}  J. Algebra 298, 460-467  (2006)

\item {[17]} M. Mecchia, B. Zimmermann, {\it On finite simple and nonsolvable groups
acting on homology 4-spheres.} Top. Appl. 153,  2933-2942  (2006)

\item {[18]} R.J. Milgram, {\it Evaluating the Swan finiteness
obstruction for finite groups.} Algebraic and Geometric Topology.
Lecture Notes in Math. 1126 (Springer 1985), 127-158

\item {[19]} J. Milnor, {\it Groups which act on $S^n$ without fixed
points.} Amer. J. Math. 79, 623-630  (1957)

\item {[20]} W. Pardon, {\it Mod 2 semi-characteristics and the converse of
a theorem of Milnor.}  Math. Z. 171,   247-268  (1980)

\item {[21]} P.A. Smith, {\it New results and old problems in finite transformation
groups.}  Bull. Amer. Math. Soc. 66, 401 - 415 (1960)

\item {[22]} M. Suzuki, {\it Group Theory II.}  Springer-Verlag 1982

\item {[23]} R.G. Swan, {\it The $p$-period of a finite group.}
Illinois J. Math. 4, 341-346 (1960)

\item {[24]} B. Zimmermann, {\it On finite simple groups acting on homology
3-spheres.}  Top. Appl. 125, 199-202 (2002)

\item {[25]} B. Zimmermann, {\it Some results and conjectures on finite
groups acting on homology spheres.} Sib. Electron. Math. Rep. 2,
233-238 (2005)   (http://semr.math.nsc.ru)

\item {[26]} B. Zimmermann, {\it Cyclic branched coverings and homology
3-spheres with large group actions.}  Fund. Math. 184, 343-353  (2004)

\item {[27]} B. Zimmermann, {\it On the classification of finite groups
acting on homology 3-spheres.} Pacific J. Math. 217, 387-395 (2004)

\bye